\documentclass[shortAfour,sagev]{sagej}

\usepackage{moreverb,url}
\usepackage{algorithm2e}

\usepackage{lipsum}
\usepackage{mathtools}

\usepackage{cuted}

\usepackage[colorlinks,bookmarksopen,bookmarksnumbered,citecolor=red,urlcolor=red]{hyperref}

\newcommand\BibTeX{{\rmfamily B\kern-.05em \textsc{i\kern-.025em b}\kern-.08em
T\kern-.1667em\lower.7ex\hbox{E}\kern-.125emX}}

\def\volumeyear{2016}

\begin{document}

\runninghead{Sizikov and Sidorov}

\title{Discrete Spectrum Reconstruction using Integral Approximation Algorithm}

\author{Valery Sizikov\affilnum{1} and Denis Sidorov\affilnum{2}}

\affiliation{\affilnum{1}ITMO University, Saint Petersburg, Russia;\\
\affilnum{2}Energy Systems Institute SB RAS, Irkutsk National Research Technical University, Irkutsk, Russia; Hunan University, Changsha, China}

\corrauth{Denis Sidorov, Energy Systems Institute of Russian Academy of Sciences,
Lermontov Street 130, Irkutsk 664033, Russia.}
\email{dsidorov@isem.irk.ru}

\begin{abstract}
An inverse problem in spectroscopy is considered. The objective is to restore the discrete spectrum from observed spectrum data, taking into account the spectrometer's line spread function. The problem 
is reduced to solution of a system of linear-nonlinear equations (SLNE)  with respect to 
intensities and frequencies of the discrete spectral lines. The SLNE is linear with respect to lines' intensities and nonlinear with respect to the lines' frequencies. The  integral approximation algorithm is proposed for the solution of this SLNE. The algorithm combines solution of linear integral equations 
with solution of a system of linear algebraic equations and avoids nonlinear equations. Numerical examples of the application of the technique,  both to synthetic and experimental    spectra, demonstrate the efficacy of the proposed approach in enabling an effective enhancement of the spectrometer's resolution.

\end{abstract}


\maketitle

\section{Introduction}

Spectral analysis is widely used for the qualitative and quantitative study of substances\cite{l1,l2,l3,l4,l5,l6}. The spectrum $u(\nu)$ characterizes dependency of  intensity of radiation as a function of frequency $\nu$. There are different types of spectra, including continuous, discrete, etc. playing principal role in scattering theory. There are various spectroscopic tools available for light decomposition into a spectrum\cite{l2}. 

Two possible avenues are available to increase the spectrometer's resolution.
The first avenue mainly concerns hardware, whereby engineers design a more sophisticated and expensive spectrometer. This paper follows the second avenue and employs mathematical processing of 
the measurements.

The purpose of this work is to propose a novel approach capable of restoring the discrete spectrum based on measured, possibly smoothed and noisy, spectral data and taking into account the  spectrometer's known line spread function. The problem is reduced 
to solution of a system of linear and nonlinear equations (some equations are linear in their variables, while some are not)  using the integral approximation algorithm. The algorithm is implemented in Matlab. A discrete (or line) spectrum  consists of discrete nearly monochromatic lines. 
The discrete spectrum is widely used to characterize  diffuse interstellar nebulae, 
low temperature gas-discharge plasma and any substance  in a deep vacuum.

Spectrum recovery (or deconvolution) is a well-known inverse problem in spectroscopy
and belongs to the ill-posed category of problems.
It is one of the  variants of the 
classic Rayleigh reduction problem~\cite{l1,l6}. 

The conventional approach to attack this challenging
inverse problem involves the Fourier self-deconvolution (FSD) technique~\cite{l31} to enhance the
resolution of spectral overlapping lines. The main idea of FSD is to employ the Fourier-transformed spectroscopy when the interferogram spectrum is measured  followed by the Fourier transform. In order to enhance the resolution FSD uses apodization (truncation of the interferogram). The objective of FSD is overlapping lines separation, i.e. when the low resolution caused by lines proximity and width, but it is not caused by  instrumental function's width as in the present paper.  It is to be noted that apodization causes lines widths artificial reduction, i.e. true lines profiles are distorted (to improve their resolution). In present paper (as well as in \cite{l6}, \cite{l8,l9,l10}) the slightly different problem is addressed: effective enhancement 
of spectrometer's resolution without artificial reduction of lines profiles.

\section{Problem statement}

Let $u(\nu)$ be  a spectrum measured e.g. by a Fabry-P\'erot interferometer~\cite{l7}.
The measured spectrum $u(\nu)$ is usually different from  true spectrum  $z(\nu)$ due to  measurement errors\cite{l6,l9} and the influence of the spectrometer's instrumental function\cite{l1,l5,l6,l8,l9}.
The following definition of the instrumental function\cite{l1,l5,l6,l8,l9} can be formulated.\\
\textit{Instrumental response function} (IF, or spectral sensitivity, transmission function, point spead function, frequency response)
$K(\nu,\nu')$ is the spectrometer's response (in terms of measured intensity) on a discrete line of unit intensity and frequency $\nu'$ when the spectrometer is tuned to the desired
frequency
 $\nu$
(proportional to a wavenumber).
Fixing $\nu$ and changing $\nu'$, the dependency $K(\nu,\nu')$ 
can be obtained in the form of a curve as shown in Figure 1.
\begin{figure}[htbp]
\centering
\includegraphics[scale=0.6]{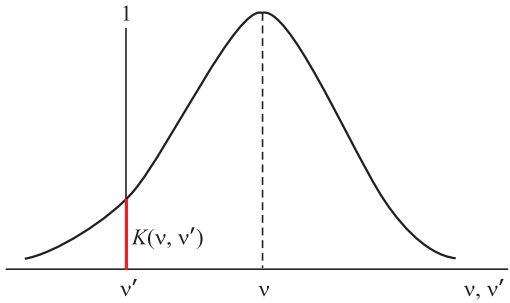}
\caption{Dependency $K(\nu,\nu')$ with concrete $\nu.$}
\end{figure}
 Similar curves can be obtained for other values of $\nu$. The result is a two-dimensional function $K(\nu,\nu')$. A wider instrumental function results in a smoother measured spectrum $u(\nu)$ when compared with the actual spectrum $z(\nu).$

The discrete spectrum $z(\nu)$  is composed of individual substantially monochromatic lines characterized by their frequencies and intensities (amplitudes). The problem of restoring 
the true discrete spectrum can be described by the following relation\cite{l6,l10,l11}
\begin{equation}
{\mathcal A} z\equiv \sum_{j=1}^n K(\nu_i, \nu_j')z_j+F=\tilde{u}({\nu}_i), \, i=\overline{1,m},\, 
\label{eq1}
\end{equation}
where $c\leq \nu_i \leq d,$ $z_j$ is intensity (amplitude) of $j$th line, $\nu_j'$ is its frequency, $n$ is 
number of lines, $\nu_i$ is digital readout of frequency $\nu$ based on spectrometer
tuning characteristics, $m$ is the number of such readouts, $[c,\,d]$ is frequency band;
source function of Eq. \ref{eq1} $\tilde{u}(\nu_i) = u(\nu_i)+\delta u(\nu_i),$
where $\delta u$ is a random component of measurement noise; $F$ is a deterministic noise component, and ${\mathcal A}$ is linear-nonlinear operator.
%\cite{st} outlined that inversed spectroscopy problem can be reduced to the Fredholm integral equation of the first kind of convolution type.

The functions $\tilde{u}(\nu_i), \, K(\nu,\nu'), \nu_i, c, d, m$ are known in system~\ref{eq1}, and
$z_j$ are the desired values. System~\ref{eq1} contains both linear and nonlinear equations,
namely the system is linear with respect to $z_j$ and $F$ and nonlinear with respect to $\nu_j'.$

System~\ref{eq1} can be considered as a system of nonlinear equations (SNE) and the known methods\cite{l12,l13} for solution of SNEs such as the Newton-Kantorovich method,  gradient, chords, gradient projection,  to name a few, can be applied directly with or without imposing  restrictions on the desired solution.

 However, these methods do not take into account the specifics of the system of linear and nonlinear Eqs~\ref{eq1} (SLNE). As a result, these methods require high computing time and memory. Moreover,  false lines may appear as
result of nonlinear system roots. Finally, the number of true spectral lines $n$ remains unknown.

The alternative avenue is to employ the folowing methods for solution of the SLNE:
\begin{enumerate}
\item[(i)] Prony method\cite{l14}. It is suitable  for SLNE with Vandermonde matrix only,
that is when matrix $K(\nu,\nu')$ varies along its lines with geometric progression. 
That is not the case for the matrix of system~\ref{eq1}.

\item[(ii)]  Peebles--Berkowitz algorithm\cite{l15}. It approximates an IF with a Taylor series resulting in errors in the solutions.
\item[(iii)] Falkovich--Konovalov algorithm\cite{l16} which is a cumbersome algorithm, difficult to implement and tune.
\item[(iv)] Golub-Mullen-Hegland method\cite{l17,l18,l19}  of variable projection. In this algorithm the Gauss-Newton method must be used to find the desired frequencies.

\end{enumerate}

In \cite{l19},  SLNE~\ref{eq1} is solved as follows. First, 
some frequencies $\nu_j'$ and $\nu_i$ are given, the amplitudes $z_j$ are determined by the Moore--Penrose pseudoinverse $z={\mathcal A}^+\tilde{u}.$
Then, the frequencies $\nu_j'$ are approximated (refined) via solving
the nonlinear least squares problem by the Gauss--Newton or Levenberg--Marquardt methods.
Thus, the problem is solved as a linear--nonlinear system.

\section{Integral approximation algorithm} 

In order to design an efficient algorithm for the solution of  system~\ref{eq1} which
takes into account its specifics, the {\it integral approximation  algorithm}\cite{l6,l10,l11} is employed in the present paper. The  integral approximation algorithm has demonstrated its efficiency in number of signal processing
problems\cite{l10,l11} as well as in spectroscopy\cite{l6}. 
 The essence of the algorithm is as follows.

The conventional approach for the solution of integral equations is to  
use the method of quadrature by  integral  approximation by the finite sum, i.e. 
the original continuous problem being reduced to a discrete problem. Here the inverse operation
is suggested: the discrete problem~\ref{eq1} is replaced with an integral equation (IE) to be solved 
using the quadrature method on finer mesh sizes. This enables the estimation of frequencies 
using only linear operations with an accuracy that depends on using a small discretization step. 

This algorithm is performed in four steps as follows.

\textsf{1st step.} Instead of system~\ref{eq1} the following linear Fredholm 
equation of the first kind\cite{l6,l8,l9,l10,l11} 
\begin{equation}
Az \equiv \int_a^b \tilde{K}(\nu, \nu')z(\nu')\, d\nu' = \tilde{u}(\nu), \,\,\, c\leq \nu \leq d
\label{eq2}
\end{equation}
is considered with respect to function $z(\nu'),$ where $A$ is a linear integral operator.
It is assumed that instead of an exact source $u$ and kernel $K$ an approximate $\tilde{u},$ $\widetilde{K}$
are given such as $||\tilde{u} - u||\leq \delta,$ $||\widetilde{K}-K ||\leq \xi,$
where $\delta$ and $\xi$ are upper bounds of the errors in the source function $u$ and kernel $K$.
The following statement\cite{l11} is valid for SLNE~\ref{eq1} and integral Eq.~\ref{eq2}.

If $z(\nu')$ in Eq.~\ref{eq2} is a generalized function consisting of sum of the Dirac delta-functions \cite{l20}
$z(\nu') = \sum_{j=1}^n z_j \delta(\nu'-\nu_j)$ and $a \rightarrow -\infty,$ $b \rightarrow \infty,$ then Eq.~\ref{eq2} is reduced to system~\ref{eq1} with $F\equiv 0,$
i.e. there is a mutual transition between SLNE~\ref{eq1} and   Eq.~\ref{eq2}.

\textsf{2nd step.} 
The solution of Eq.~\ref{eq2} is a known {\it ill-posed problem}\cite{l6,l20, l21,l22,l23,l24,l25,l26,l27}. It is to be noted that direct application of the quadrature method results in a so-called ``saw'', thus rendering the method very unstable\cite{l6,l23}. 
Therefore, stable methods such as zero-order Tikhonov regularization method\cite{l6,l20,l21,l22,l23,l24,l25,l26,l27} should be applied. Using the Tikhonov  method we replace  Eq.~\ref{eq2} with the following regularized equation\\
%$$\alpha z_{\alpha}(\varphi) +$$
\begin{equation}
\alpha z_{\alpha}(\varphi) +\int_a^b R(\varphi,\nu')z_{\alpha}(\nu')\, d\nu' = U(\varphi), \label{eq3}
\end{equation}
where $a\leq \varphi\leq b,$ $\alpha > 0 $ is a regularization parameter, 
\begin{equation}
R(\varphi,\nu') = R(\nu', \varphi) = \int_c^d K(\nu,\varphi) K(\nu, \nu') \, d\nu,
\label{eq4}
\end{equation}
\begin{equation}
U(\varphi) = \int_c^d K(\nu, \varphi) \tilde{u}(\nu) \, d\nu.
\label{eq5}
\end{equation}

%Integral  Eq.~\ref{eq3} is supposed to be soved using quadrature method 
%with smaller regularization parameter $\alpha$ and a fine sampling step 
%$h=\Delta \nu'.$ Numeric experiments have demonstrated that 
%regularization parameter $\alpha$ should be assumed by 2-3 orders of magnitude smaller 
%comparing with conventional methods such as generalized discrepancy principle\cite{l21,l22,l26}.
% Otherwise solution will be over-smoothed. The sampling step $h$ should be as small as possible to make  several hundred discrete samples for $\nu'$. The small sampling step $h$ is important in this algorithm.

For solution of integral Eq.~\ref{eq2} quadrature methods \cite{l6,l21,l22,l23,l26} can be 
applied by replacement of the integral with a finite sum leading to the system of linear algebraic equations
(SLAE)
${ A} { z}=\tilde{ u}.$ In that case quadrature coefficients are assumed to be equal one
in SLAE $ { A} { z} =\tilde{ u}$ to make it structurally close to SLNE \ref{eq1}. 
% Paul L: I didn't understand the intended meaning of the previous sentence. The grammar is not right, but I couldn't correct it.
 In order to obtain the stable solution
instead of SLAE ${ A} { z}=\tilde{ u}$ we solve following stabilized SLAE which follows from Eqs.~\ref{eq3}-\ref{eq5}
\begin{equation}
(\alpha {   E} +{   A}^T {   A}){   z}_{\alpha} = {   A}^T \tilde{   u}
\label{eq6}
\end{equation}
with the square  matrix ${   R}=(\alpha {   E}+{   A}^T{   A})$ of dimention $N \times N,$ where
$N$ is the number of discrete samples for $\nu',$ ${   E}$ is the identity matrix, ${   A}$ is a $m \times N$
matrix based on kernel $K;$ ${   A}^T$ is the transposed matrix; $\tilde{   u}$ is a given 
vector of length $m$ which is the measured spectrum. The desired solution of SLAE \ref{eq6} is 
intensity vector ${   z}_{\alpha j} \equiv  {   z}_{\alpha}(\nu_j^{\prime}), j=\overline{1,N}$
on a fine mesh with reduced step $h=\Delta\nu' = (b-a)/(N-1). $

A critical matter is the strategy for the choice of the regularization parameter $\alpha$. The  conventional method
for choosing the parameter $\alpha$ is based on the Morozov discrepancy principle \cite{l28}.
In this principle the parameter $\alpha$ can be chosen using equality (see Fig. 3 below)
\begin{equation}
||{   A} {   z}_{\alpha}  - \tilde{   u}||=\delta,
\label{eq7}
\end{equation}
where the  discrepancy $||{   A} {   z}_{\alpha}  - \tilde{   u}||$ can be written as follows
\begin{equation}
||{   A} {   z}_{\alpha}  - \tilde{   u}||_2 := \biggl\{ \sum_{i=1}^m \biggl[\sum_{j=1}^N A_{ij}z_{\alpha j}-\tilde{u}_j  \biggr]^2 \biggr\}^{1/2}.
\label{eq8}
\end{equation}
In Eq.~\ref{eq7} $\delta = ||\tilde{   u}-{   u}||_2,$ where $\tilde{   u}$ is the measured noisy spectrum,
and $u$ is the unknown non-noisy spectrum. Here $u$ is approximated by a smoothing spline
through the noisy samples $\tilde{u}.$ Such an approach has been successfully employed
in \cite{l25,l29}. $\alpha_d$ denotes a regularization parameter chosen according to 
the discrepancy principle.

The discretization step $h$ should be as small as possible in order to make  several hundred discrete samples $N$ for $\nu',$  i.e. $N\gg n,$ where $n$ is expected number of lines (see Eq. \ref{eq1}).
The number $m$ for $\nu$ is the same as in Eq.~\ref{eq1}, or it can be  increased using the spline function. Employing a small discretization step $h$ is important in the application of this algorithm.

It is to be noted that $N \lessgtr m$, i.e. at any ratio of $N$ and $m$ the Tikhonov regularization method gives the solution (normal pseudo-solution) of the integral Eq.~\ref{eq3} and the corresponding
SLAE. If $N$ is assumed to be greater than $m$ then, based on a limited number of measurements 
$m$ (e.g., $ m \approx 100$), one can obtain a solution ${   z}_{\alpha}(\nu')$ defined on a greater 
number of nodes $N$ (e.g., $N \approx 400$) and using a smaller step $h=\Delta \nu'.$

The desired solution ${   z}_{\alpha}(\nu')$ is thus obtained. Such a solution can help to distinguish
the additional lines in the desired spectrum. False line-peaks may arise and must be filtered, and this issue is addressed below. 

The following estimate\cite{l6,l9,l10,l11} of the error's norm of the regularized solution ${   z}_{\alpha}(\nu')$ of Eq. \ref{eq3} is
 carried out
\begin{equation}
\varepsilon(\alpha) \equiv || \Delta {   z}_{\alpha}||\leq \left( \frac{||{   A}||}{2\sqrt{\alpha}}\eta + \frac{p\alpha}{p\alpha+1} \right)||{   z}_{\alpha}||,
\label{eq9}
\end{equation} 
where  $\eta = \delta_{\text{rel}}+\xi_{\text{rel}},$ $\delta_{\text{rel}}=\delta/||{   u}||,$ $\xi_{\text{rel}}=\xi/||{    A}|| $ are relative errors of input data,
and parameter $p$ can be chosen using the model spectra processing method\cite{l9}.

Let us estimate the error $\sigma_{\nu_j'}$ of  $j$-th line frequency $\nu_j'.$  The solution ${    z}_{\alpha}(\nu')$ is expanded in a Taylor series in the neighborhood of frequency $\tilde{\nu}_j'$  as follows
\begin{equation}
{   z}_{\alpha}(\nu') \approx {   z}_{\alpha}(\tilde{\nu}_j')+\frac{1}{2}{   z}_{\alpha}''(\tilde{\nu}_j')\Delta_j^2,
\label{eq10}
\end{equation}
where ${z}_{\alpha}''(\tilde{\nu}_j')=\frac{\partial^2 {   z}_{\alpha}(\nu')}{\partial \nu'^2}\bigl|_{\nu'=\tilde{\nu}_j'},$
$\Delta_j=\nu'-\tilde{\nu}_j'=\sigma_{\nu_j'}.$
It is to be noted that ${   z}_{\alpha}'(\tilde{\nu}_j')=0$ and ${   z}_{\alpha}^{'''}(\tilde{\nu}_j')\approx 0.$ That is the reason why series~\ref{eq10} should provide
an accurate approximation of ${    z}_{\alpha}(\nu')$ in the  $j$th line neighborhood.
Here $\tilde{\nu}_j^{\prime}$ is the approximate frequency of the $j$th line corresponding to 
the some peak in the solution ${   z}_{\alpha}(\nu'),$ $\nu'$ is some frequency value in the neighborhood of $\tilde{\nu}_j', $ in particular, exact frequency $\nu_j'$ value.
From Eq.~\ref{eq10} as follows:
\begin{equation}
(\nu' - \tilde{\nu}_j')^2 \approx \frac{2[{   z}_{\alpha}(\nu') - {   z}_{\alpha}(\tilde{\nu}_j') ]}
{{   z}_{\alpha}''(\tilde{\nu}_j')} 
\label{eq11}
\end{equation}  
or using absolute values
\begin{equation}
|\nu' - \tilde{\nu}_j'| \approx  \sqrt{\frac{2|{   z}_{\alpha}(\nu') - {   z}_{\alpha}(\tilde{\nu}_j') |}
{|{   z}_{\alpha}''(\tilde{\nu}_j')|}}. 
\label{eq12}
\end{equation}  
However, such an error estimate of the $j$th line frequency can not be employed because 
$\nu',$ is not known, in particular, the exact value of frequency $\nu_j'.$

A more constructive approach for frequency error estimation is to employ the norms.
Let us make the notation $||{   z}_{\alpha}(\nu') - {   z}_{\alpha}(\tilde{\nu}_j')|| = 
||\Delta {   z}_{\alpha}|| = \varepsilon(\alpha)$ and obtain the following estimate 
of the $j$th line frequency $\nu_j'$ error  using $\varepsilon(\alpha)$ estimate 
 the norm of the regularized solutions
\begin{equation}
\sigma_{\nu_j'} \approx \sqrt{\frac{2\varepsilon(\alpha)}{|{   z}_{\alpha}''(\tilde{\nu}_j')|}},
\label{eq13}
\end{equation}
where $\varepsilon(\alpha)$ is calculated using Eq.~\ref{eq9}. 

%Using notation $||z_{\alpha}(\nu') - z_{\alpha}(\nu_j')|| = ||\Delta z_{\alpha} || = \sigma(\alpha) $ the following 
%estimate of error of frequency $\nu_j'$ of $j$th line can be obtained
%\begin{equation}
%\sigma_{\nu_j'} \approx \sqrt{\frac{2\sigma(\alpha)}{|z_{\alpha}''(\nu_j')|}},
%\label{eq14}
%\end{equation}
%where $\sigma(\alpha)$ is selected according to Eq.~\ref{eq6}.

Error $\sigma_{\nu_j'}$ also grows due to the finite discretization step $h=\Delta\nu'$ used for numerical solution of the integral Eq.~\ref{eq2}. The error growth is approximately equal to $h/2.$
Finally, taking the square root of two error terms the following more accurate error estimate can be obtained

\begin{equation}
\sigma_{\nu_j'} \approx \sqrt{\frac{2\varepsilon(\alpha)}{|z_{\alpha}''(\tilde{\nu}_j')|}+\left(\frac{h}{2}\right)^2}.
\label{eq14}
\end{equation}
 Eq.~\ref{eq14} ${   z}_{\alpha}''(\tilde{\nu}_j')$ is an estimate of the second derivative of the regularized 
solution ${   z}_{\alpha}(\nu')$ for $\nu'=\tilde{\nu}_j',$ i.e. in the area of $j$th line. If $h={\text{const}},$
then 
\begin{equation}
{   z}_{\alpha}''(\tilde{\nu}_j') \approx \frac{{   z}_{\alpha}(\tilde{\nu}_j' -h) - 2{   z}_{\alpha}(\tilde{\nu}_j')+{   z}_{\alpha}(\tilde{\nu}_j'+h)}{h^2}.
\label{eq15}
\end{equation}

Estimate~\ref{eq14} shows that the determination error of frequencies $\sigma_{\nu_j'}$  decreases as the  discretization step $h=\Delta \nu'$ decreases  and as the 
second derivative $z_{\alpha}''(\tilde{\nu}_j')$ increases, which increases with decreasing values of the regularization parameter $\alpha$. That is the reason why IE~\ref{eq3} should be solved using  a
small  step $h$ and a small regularization parameter $\alpha.$

The following error metrics can be employed for analysis of 
the model spectrum with known exact intensities $z_j$ and frequencies $\nu_j'$
$(j=\overline{1,n}).$
\begin{enumerate}
\item
Root mean square error (RMSE) of lines' intensities (cf. \cite{l27}):
\begin{equation}
\varepsilon = \biggl \{ \frac{1}{n} \sum_{j=1}^n (\tilde{z}_j - z_j)^2  \biggr \}^{1/2}  = ||\tilde{   z} - {   z}||_2,
\label{eq16}
\end{equation}
\item RMSE of lines' frequences:
\begin{equation}
\xi = \biggl \{ \frac{1}{n} \sum_{j=1}^n (\tilde{\nu'}_j - \nu_j')^2  \biggr \}^{1/2}   = ||\tilde{   \nu}' - \nu'||_2.
\label{eq17}
\end{equation}
\end{enumerate}
RMSE metrics depend on the system of units making them difficult to compare. 
The following relative  root mean square errors (RRMSE) are used in our experiments:
\begin{equation}
\varepsilon_{\text{rel}} = \biggl \{ \frac{ \sum_{j=1}^n (\tilde{z}_j - z_j)^2}{\sum_{j=1}^n z_j^2}  \biggr \}^{1/2}  = \frac{||\tilde{   z} - {   z}||_2}{||z||_2},
\label{eq18}
\end{equation}
\begin{equation}
\xi_{\text{rel}} = \biggl \{\frac{  \sum_{j=1}^n (\tilde{\nu'}_j - \nu_j')^2}{\sum_{j=1}^n \nu_j'^2}  \biggr \}^{1/2}   = \frac{||\tilde{\nu}' - \nu'||_2}{||\nu'||_2},
\label{eq19}
\end{equation}
where $\tilde{z}_j,$ $\tilde{\nu}_j'$ are approximated intensities and frequencies, and
${z}_j,$ ${\nu}_j'$ are exact values of intensities and frequencies. 
Total RRMSE can be computed as follows:
\begin{equation}
\zeta_{\text{rel}}  = \sqrt{\varepsilon_{\text{rel}}^2 +\xi_{\text{rel}}^2 }.
\label{eq20}
\end{equation}
In the examples below the  RMSEs $\varepsilon_{\text{rel}},$
$\xi_{\text{rel}},$ $\zeta_{\text{rel}}$ are employed.

\textsf{3rd step}. $L $ most powerful lines (maxima) are selected 
in the resulting solution $z_{\alpha}(\nu')$ on the basis of additional information.
Here $L$ is given with a margin such as $L \geq n,$ but $L\ll N,$ where $n$ is the expected number of lines.
The frequencies of the most powerful maxima $\tilde{\nu}_j', \, j=\overline{1,L}$ are recorded.

\textsf{4th step}. Solution of the following refined system of linear algebraic equations (SLAE) of lower
order
\begin{equation}
\sum\limits_{j=1}^L K(\nu_i, \tilde{\nu}_j')\tilde{z}_j + \widetilde{F} = \tilde{u}(\nu_i),\,
i=\overline{1,m}, c\leq \nu_i\leq d.
\label{eq21}
\end{equation}
This system has more equations than unknowns $\tilde{z}_j,$ i.e. system~\ref{eq21} is overdetermined since $ m >L.$
Usually, $m \sim 10^2, L\sim 10^1$ and the least squares method (LSM) have to be used for such SLAE without regularisation with respect to line intensities $\tilde{z}_j$ and  $\tilde{F}$
taking into account $\tilde{\nu}_j'$ defined in Step 2 and 3.
In the LSM the SLAE with a $(L+1) \times (L+1)$ square matrix is obtained.
$\tilde{z}_j$ and $\tilde{F}$ are supposed to be chosen according to the following threshold 
\begin{equation}
\tilde{z}_j \geq Z, \, j=\overline{1,k}, \tilde{F}>0,
\label{eq22}
\end{equation}
where $Z>0$ is the threshold value, and $k\leq L$ is the number of $\tilde{z}_j$ passing the 
threshold $Z.$ 
The threshold value can be chosen  as follows~\cite{l10,l11}
\begin{equation}
Z = \delta \sqrt{-2\ln P_{\text{fa}}},
\label{eq23}
\end{equation}
where $P_{\text{fa}} \in [0,\,1]$ is given conditional probability of a false alarm. However, 
false maxima  are negative or nearly zero, so there is no need to use  Eqs.~\ref{eq22}, \ref{eq23}.

For sake of clarity let us summarize the algorithm.

\begin{algorithm} {\bf Stage 1.} IE.~\ref{eq2} is considered instead of 
SLNE~\ref{eq1}.\\
				    {\bf Input:} $\tilde{u}(\nu), n, m, c, d, K(\nu,\nu'), a, b, h,\delta, N$.\\
				%	$n_{min}$ -- ìèíèìàëüíîå êîëè÷åñòâî òî÷åê â âåòâè; \\
				{\bf Stage 2.} IE~\ref{eq2} is solved using Tikhonov regularization (Eqs. 
\ref{eq3} -- \ref{eq5}) using reduced discretization step $h$ to get the SLAE (of $m$ equations) with respect to 
$N\gg n$ desired $z_{\alpha j} \equiv z_{\alpha}(\nu_j'),$ where regularization parameter 
$\alpha$ is chosen using the Morozov discrepancy principle \ref{eq7} -- \ref{eq8}.
The error of regularized solution is computed.\\
{\bf Output:} $z_{\alpha}(\nu').$\\
{\bf Stage 3.} Define $L$ maxima peaks in  $z_{\alpha}(\nu'),$ $L\geq n$ and 
$L\ll N.$\\
{\bf Output:} {\it frequencies $\tilde{\nu}_j', j=\overline{1,L}.$}\\
{\bf Stage 4.} Refined SLAE \ref{eq21} of lower order $m\times (L+1) $ is solved to find
lines intensities $\tilde{z}_j$ and background $\tilde{F}$ using Eq.~\ref{eq22}.\\
{\bf Output:} {\it lines intensities $\tilde{z}_j$, background $\tilde{F}$ and lines number $k$.}

\end{algorithm}

The advantage of the algorithm is that it allows linearization of the nonlinear problem of determination of the frequencies of spectral lines $\tilde{\nu}_j'$ using the linear IE~\ref{eq2}.

\section{Numerical experiment}

The proposed algorithm for discrete spectrum reconstruction is implemented as a software package in Matlab. In this section the package is tested on both
  model (for testing and verification) spectra and real spectra.

\subsection{1st example: model spectrum}

Following  paper\cite{l9},  IF $K$ is assumed to have a variable width depending on frequency $\nu,$ i.e. in general the IF is not a convolution type kernel and $K = K(\nu,\nu').$ 
Such IFs are typical in broadband spectroscopy.
The function $\tau(\nu)$, the half-width of the IF on half-power, is conventionally used. It is proportional 
to wavelength  $\lambda$ and inversely proportional to the frequency, i.e. $\tau(\nu)=q/\nu,$ 
where the coefficient $q$ numerically determines  the half-width of the IF and depends on the particular spectrometer design and setup.

The following most typical IFs\cite{l1,l9,l19} have been considered. These IFs are
characterized by their half-widths $\tau(\nu)$ due to diffraction, aberrations and slots types.
%\begin{enumerate}

\noindent 1. Slotlike (rectangular) IF characterized by wide slot only
\begin{equation}
  K(\nu,\nu')=\left\{
            \begin{array}{ll}
              \frac{1}{2\tau(\nu)}, & \hbox{if} |\nu - \nu'|\leq \tau(\nu), \\
              0, & \hbox{otherwise}.
            \end{array}
          \right.
\label{eq24}
\end{equation} 

\noindent 2. Triangular IF (also takes into account  slit distortions only, it is the convolution of two rectangular IFs)

\begin{equation}
 K(\nu,\nu')=\left\{
            \begin{array}{l}
              \frac{1}{2\tau(\nu)} \left(1-\frac{|\nu-\nu'|}{2\tau(\nu)}  \right ), \,\hbox{if\,} |\nu - \nu'|\leq 2\tau(\nu), \\
              0, \,\hbox{otherwise.}
            \end{array}
          \right.
\label{eq25}
\end{equation} 

\noindent 3. Rayleigh diffraction IF (in this case  device has an infinitely narrow slot and  the IF 
is determined by diffraction at a rectangular aperture diaphragm)
\begin{equation}
   K(\nu,\nu')=\frac{1}{\gamma(\nu)} \left\{\frac{\sin[\pi(\nu-\nu')/\gamma(\nu)]}{\pi(\nu-\nu')/\gamma(\nu)}  \right\}^2,
\label{eq26}
\end{equation}
where $\gamma(\nu)=2\tau(\nu)/0.8859=2.2576\tau(\nu).$

\noindent 4. Gaussian IF\cite{l6,l20,l30} (monochromator case where the diffraction and aberration distortions are not very large)
\begin{equation} 
 K(\nu,\nu')=\frac{1}{\sqrt{2\pi}\sigma(\nu)} \exp \left( - \frac{(\nu - \nu')^2}{2\sigma^2(\nu)} \right ),
\label{eq12add}
\end{equation}
where $\sigma(\nu) = \tau(\nu)/\sqrt{2\ln 2} = 0.8494\tau(\nu)$\cite{l6}.

\noindent 5. Dispersion or Lorenz IF (for spectrographs of  small, medium and sometimes large dispersion)
\begin{equation} 
   K(\nu,\nu')= \frac{\tau(\nu)/\pi}{(\nu-\nu')^2+\tau^2(\nu)}.
\label{eqDT}
\end{equation}

\noindent 6. Exponential IF (of photosensitive layer)
\begin{equation}
K(\nu,\nu') = \frac{\ln 2}{2\tau(\nu)}\exp\left(-\frac{\ln 2}{\tau(\nu)}|\nu - \nu'| \right ).
\label{eq13add}
\end{equation}

\noindent 7. Voigt IF (convolution of Gauss and Dispersion IFs)\cite{l32}. Its explicit formula is omitted for the sake of brevity.\\

\noindent Remark.\\
Obviously IFs can also be written in terms of wavelength $\lambda$ (wavenumber $\nu/c = 1/\lambda$), e.g. Eq.~\ref{eq13add} can be written\cite{l9} as 
\begin{equation}
K(\lambda,\lambda') = \frac{\ln 2}{2\tau(\lambda)}\exp\left(-\frac{\ln 2}{\tau(\lambda)}|\lambda - \lambda'| \right ),
\label{eq30}
\end{equation}
where $\tau(\lambda) = q\lambda$ is half-width of IF on half-power.

It is useful to employ the half-widths $\tau(\nu)$ when it comes to comparison of IFs with 
the same $\tau(\nu).$ Such a comparison was conducted\cite{l9} for continuous spectra.
It is demonstrated that the most accurate results can be achieved for 
dispersion and exponential IFs and the least accurate results were achieved for 
slot and triangular IFs. This observation arose from processing the discrete spectra with various IFs \ref{eq24}-\ref{eq13add}
using the integral approximation algorithm. 

 Application of the integral approximation algorithm to the discrete spectrum described above confirmed this observation.

In order to demonstrate the efficiency of the proposed approach the following model example\cite{l6,l10,l11} with synthetic data is considered in this section.

The true spectrum is given as seven discrete spectral lines with the following amplitudes (in arbitrary units, a.u.) $z_1=4.4,$ $z_2 = 4.6,$ $z_3 = 1.1,$ $z_4=3.2,$
$z_5=3.2,$ $z_6=2.8,$ $z_7=3.6$ and frequencies (given in a.u. as well)
$\nu_1'=2.28,$ $\nu_2'=2.36,$ $\nu_3'=2.95,$ $\nu_4'=3.02,$ $\nu_5'=3.56,$
$\nu_6'=3.64,$ $\nu_7'=3.69.$
The spectrometer's IF is given\cite{l9} as the following frequency--non-invariant
Gaussian function (by Eq.~\ref{eq12add}), in which width decreases with increasing frequency $\nu$ 
\begin{equation}
K(\nu, \nu') =  \frac{g}{\sqrt{2\pi}\sigma(\nu)} \exp \left(-\frac{(\nu-\nu')^2}{2\sigma^2(\nu)} \right),
\label{eq31}
\end{equation}
where $\sigma(\nu)=\sigma_0\sqrt{1-0.16 \nu}, \sigma_0=0.05, $ $ g=0.075$ is a normalizing factor, all  in a.u. For each kind of spectrometer $\sigma(\nu), \, \sigma_0$ and $g$ will be specific.
The level of the deterministic noise $F$ was assumed to be
equal $0.2$, and the random measurement noise  has a standard deviation
 (SD) equal to $0.05$ ($\approx 2\%$ noise). 

The limits were selected as follows
$a=c=2,$ $b=d=4$ a.u. (see Eqs. \ref{eq2}-\ref{eq5}), $m=101$ is the number of 
``experimental'' samples on $\nu$ (see Eq. \ref{eq1}); $N=401$ is the number of samples on
$\nu'$ in the SLAE \ref{eq6} solution, discretization step $h=\Delta \nu' = (b-a)/(N-1).$

Figure 2 shows: the true discrete (line) spectrum $z(\nu)$ consisting of 7 lines, the measured (experimental) noise free spectrum $ u(\nu)$, the noisy spectrum $\tilde{u}(\nu),$
and the cubic spline approximation  $\tilde{u}(\nu)$ is realised (length $m=101$).

In addition, the IF of the spectrometer $K(\nu, \nu')$ according to  Eq. \ref{eq12} is given at low and high frequencies $\nu$. 
It can be observed that the true spectrum $z(\nu)$ has closely spaced lines (two lines coorresponding low frequency in Figure 2, two lines are in the middle and three lines corresponding high frequency)  which 
remains unrecognized in the 
 measured spectrum $u(\nu)$.

\begin{figure}[htbp]
\includegraphics[scale=0.465]{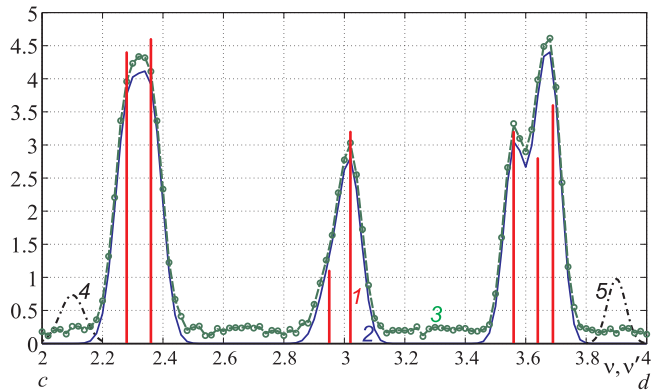}
\caption{Numerical example of the direct problem.
Frequencies $\nu$ and $\nu'$ are located on the horizontal axis, $z, u, \tilde{u},$ and $ K$ are on the vertical axis in a.u. ``{\it 1}'' is the true discrete spectrum, ``{\it 2}'' is the measured spectrum ${   u}(\nu),$
``{\it 3}'' is the measured noisy spectrum $\tilde{   u}(\nu),$ (dashed) and smoothing spline (circles)  ``{\it 4}'' is the IF $K(2.1,\nu' )$ on low frequency $\nu$, ``{\it 5}'' is $K(3.9,\nu')$ is the IF on high frequency $\nu$. }
\end{figure}

In this experiment the integral approximation algorithm is applied for the true spectrum. 
First,  IE~\ref{eq2} is supposed to be solved using the Tikhonov regularization method 
 following Eqs.~\ref{eq3}-\ref{eq5} by means of the SLAE \ref{eq6}
solution.  The regularization parameter $\alpha$ is chosen using the 
discrepancy principle \ref{eq7} - \ref{eq8} with $\alpha = \alpha_d = 10^{-1.4} = 0.04$
as shown in Fig. 3. 

\begin{figure}[htbp]
\includegraphics[scale=0.48]{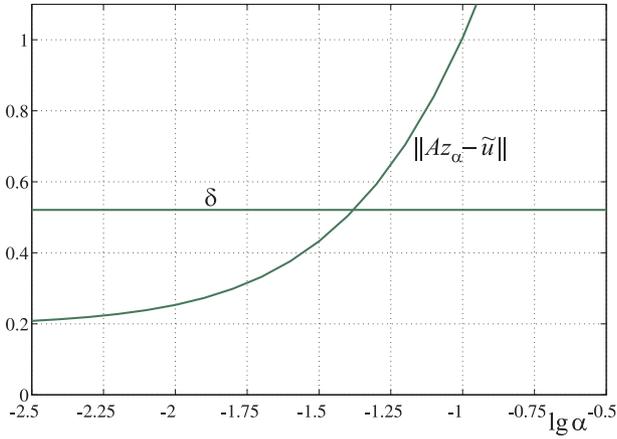}
\caption{Discrepancy principle for specifying the regularization parameter $\alpha$. }
\end{figure}

Fig. 4 shows the noisy spectrum  $\tilde{u}(\nu)$ according to Eq.~\ref{eq1} for $m=101,$
i.e. for relatively small (insufficient) samples $m.$

% Here fig 4 goes

\begin{figure}[htbp]
\includegraphics[scale=0.467]{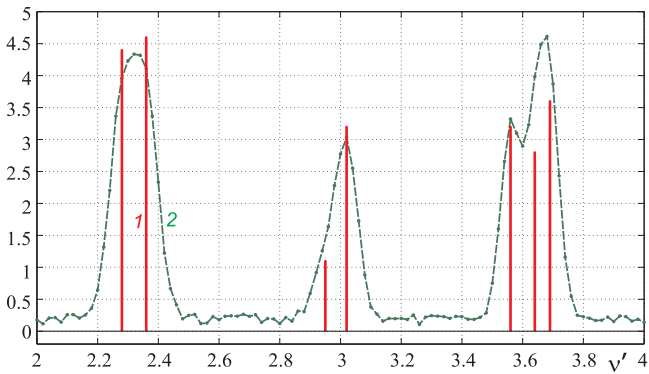}
\caption{Measured and exact spectra. ``{\it 1}'' is the exact spectrum ${   z}(\nu')$ (vertical solid lines);
``{\it 2}'' is the measured spectrum $\tilde{   u}(\nu)$ of length $m = 101.$}
\end{figure}

Using the spectrum $\tilde{   u}(\nu)$ of length $m = 101$ the regularized solution 
${   z}_{\alpha}(\nu)$ is obtained, its length $N = 401$ for $\alpha = \alpha_d$
(see Fig. 5). One may observe here  the lack of resolution
(the two rightmost high frequency lines are missing), so steps 3 and 4 of the integral approximation algorithm 
are not carried out.

% Here fig 5 goes

\begin{figure}[htbp]
\includegraphics[scale=0.46]{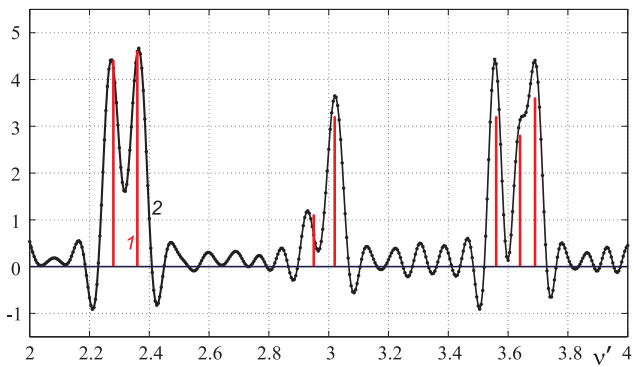}
\caption{Numerical example. Inverse problem with insufficient $m=101.$
``{\it 1}'' is the true spectrum $z(\nu')$ (vertical solid lines); ``{\it 2}'' is the regularized
solution $z_{\alpha}(\nu')$ of length $N=401.$}
\end{figure}

In order to increase the effective resolution a cubic spline ${   u}_s(\nu)$ of length $m=401$
is used to approximate the measured spectrum $\tilde{u}(\nu)$ as shown in Fig.~6.

%fig 6

\begin{figure}[htbp]
\includegraphics[scale=0.46]{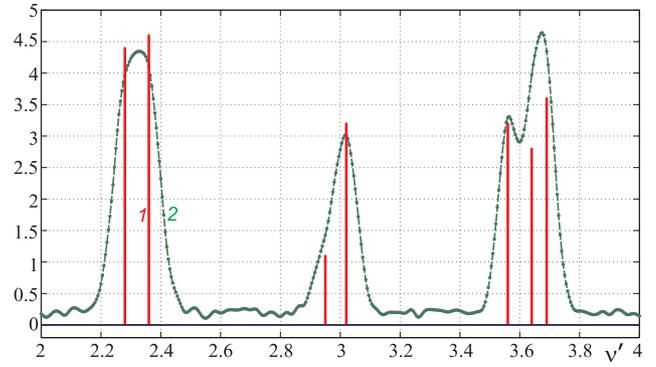}
\caption{Measured and exact spectra. Here ``{\it 1}'' is the exact spectrum ${   z}(\nu')$ (vertical solid lines); ``{\it 2}'' is a smoothing spline ${   u}_s(\nu)$ of length $m=401.$ }
\end{figure}

The spline $u_s(\nu)$ of length $m = 401$ is used to approximate 
the measured spectrum $\tilde{   u} (\nu).$ Then the regularized solution is calculated using
$\alpha = \alpha_d$ (see. Fig. 7).

It can be observed that in spite of the fact that in Fig. 7 the length of the solution is the same as in Fig. 5
($N = 401$), the resolution is improved significantly. This is due to the fact that the  length of the measured spectrum $m$ is increased by using the spline.

 %The resulting regularized solution $z_{\alpha}(\nu)$ is shown in Figure~3. All the
%true spectral lines are recognized, but there are many false lines (peaks).

In the regularized solution $z_{\alpha}(\nu)$ in Fig. 7 the first $L=12$ most powerful highs are 
taken and their frequences $\tilde{\nu}_j',$  $j = \overline{1,L}$  are recorded. 
Then the refining SLAE~\ref{eq21} is solved using LSM (without regularization) with respect to $L+1=13$ unknowns: 12 amplitudes $\tilde{z}_j$ and 
$F$. Figure 7 shows the reconstructed spectrum values $z_j(\tilde{\nu}_j')$. 
It is to be noted that  all the false spectral lines correspond to negative values or values close to zero. Here the threshold value $Z = 0.2F$ (20\%), which has not been used.
The true maximum values $\tilde{z}_j$ are indeed very close to the true amplitudes $z_j.$
The RRSME errors are $\varepsilon_{\text{rel}}=0.0621,$ $\xi_{\text{rel}} = 0.0028$
and $\zeta_{\text{rel}}=0.0622.$

\begin{figure}[htbp]
\includegraphics[scale=0.41]{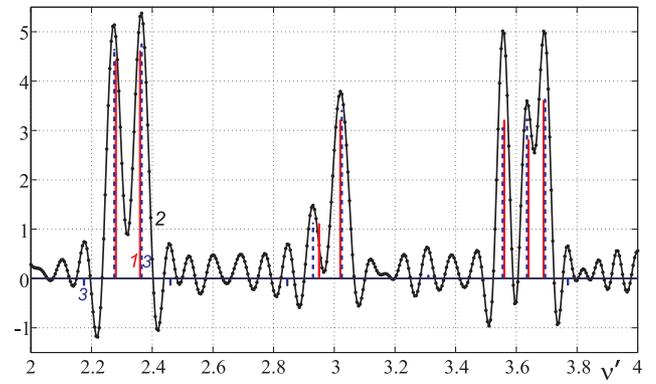}
\caption{Numerical example of the inverse problem with increased $m=401.$
The frequency $\nu'$ lies on the horizontal axis, $z$ is on the vertical axis in a.u. ``{\it 1}'' is the true discrete spectrum ${   z}(\nu')$ (vertical solid lines), ``{\it 2}'' is the regularised solution ${   z}_{\alpha}(\nu')$ of length $N=401$ (solid line);
``{\it 3}'' is the reconstructed spectrum ${   z}_{j}(\tilde{\nu}_j')$ (vertical  dashed lines). }
\end{figure}

All the 7 spectral lines have been perfectly recognized as well as their frequencies $\tilde{\nu}_j'$ and
intersities $\tilde{z}_j$ in the model example. This was achieved despite the noise and the IF which was
selected to demonstrate the efficiency of the integral approximation algorithm.

For sake of comparison let us now solve this synthetic example using   the variable projection Golub-Mullen-Hegland method\cite{l17,l18,l19}. The efficient implementation of 
variable projection method\cite{addref} (function {\tt varpro.m}) was employed. Following
notations\cite{addref}, it's assumed {\tt Ind = dPhi = []} which means that the derivatives (Jacobian matrix) were not used for solving both linear and nonlinear problems. This can reduce the accuracy of the solution, but  simplifies the solution procedure.
The problem was addressed using  three  versions as follows.

In the 1st version the  uniform mesh was used on $\nu$={\tt [2:0.1:4]}, number of nodes $m=21.$ Initial guess $\nu'$={\tt [2.26, 2.38, 2.93, 3.04, 3.54, 3.63, 3.71]}, number of lines $n=7.$
Resulting $z$ and $\nu'$ are listed in the table below. The significant error can be observed
especially in linear parameters $z$ (lines intensities). Apparently such error is caused by missing
initial guess in the function {\tt varpro.m}.

In the 2nd version the non-uniform mesh was used on $\nu$={\tt [2.1; 2.2; 2.25; 2.3; 2.35; 2.4; 2.5; 2.6; 2.8; 2.9; 2.95; 3; 3.05; 3.1; 3.2; 3.3; 3.4; 3.5; 3.55; 3.6; 3.65; 3.7; 3.75; 3.8; 3.9]}, number of nodes $m=25.$ Initial guess for $nu'$ was selected same as in the 1st version as described above. Table 1 demonstrates errors reduction and they are close to errors of  
the integral approximation algorithm: 
$\varepsilon_{\text{rel}}=0.0621,$ $\xi_{\text{rel}} = 0.0028$
and $\zeta_{\text{rel}}=0.0622.$

The true number of spectral lines is normally remains not known. That's the reason why in the 3rd version
we assumed $n=9$ (instead of true $n=7$). Initial guess $\nu'$ = {\tt [2.1, 2.26, 2.38, 2.93, 3.04, 3.1, 3.54, 3.63, 3.71]} (two additional lines were introduced). The no-nuniform
mesh was selected  same as in the 2nd version. In order to calculate errors 
$\varepsilon_{\text{rel}},$ $\xi_{\text{rel}} $
and $\zeta_{\text{rel}}$ based on Eqs. \ref{eq18} -\ref{eq20} two false lines with zero intensities
were considered:
$z$={\tt [0; 4.4; 4.6; 1.1; 3.2; 0; 3.2; 2.8; 3.6], } $\nu'$={\tt [2.1; 2.28; 2.36; 2.95; 3.02; 3.1; 3.56; 3.64; 3.69]. }
 Table 1 shows that the intensities $z$ of false lines are correspondingly equal to 0.407 and 0.767, although one would expect near-zero values. It can be concluded that the false lines suppression
does not work here.

\begin{strip}
\begin{center}
\begin{tabular}{ l|l }
\hline
  \hline
 1st version & error\\
  \hline
   $z:$ {\tt 6.413, 3.195, 1.155, 4.421, 1.467, 12.09, 0.682}& $\varepsilon_{\text{rel}}=1.1272$ \\
  $\nu':$ {\tt 2.279, 2.417, 2.863, 3.026, 3.465, 3.652, 3,847} & $\xi_{\text{rel}}=0.0257$ \\
  {\tt Warnings: linear parameters are not well-determined} & $\zeta_{\text{rel}}=1.1275$ \\
  Number of iterations $N=28,$ time $T=1.5$~sec. & $\,$\\
\hline
  \multicolumn{2}{l}{2nd version } \\
\hline
  $z:$ {\tt 4.848, 4.783, 1.361, 3.546, 3.271, 3.607, 3.498} & $\varepsilon_{\text{rel}}=0.1146$ \\
  $\nu':$ {\tt 2.280, 2.363, 2.943, 3.020, 3.554, 3.638, 3.696} &  $\xi_{\text{rel}}=0.0014 $ \\
  {\tt Warnings: there are no.} & $\zeta_{\text{rel}}=0.1146 $ \\
  $N=5,$ $T=0.3$ sec. & $\,$ \\
  \hline
  \multicolumn{2}{l}{3rd version } \\
\hline
$z:$ {\tt 0.407, 4.859, 4.743, 1.302, 3.555, 0.767, 3.271, 3.607, 3.498} & $\varepsilon_{\text{rel}}=0.1481$\\
$\nu':$ {\tt 2.141, 2.280, 2.363, 2.941, 3.019, 3.145, 3.554, 3.638, 3.696} &
$\xi_{\text{rel}}=0.0069$\\
  {\tt Warnings: there are no.} & $\zeta_{\text{rel}}=0.1483 $ \\
  $N=7,$ $T=0.4$ sec. & $\,$ \\
  \hline
\hline
\end{tabular}
\vspace*{0.3cm}

{\bf Table 1.} Results of the variable projection method of Golub-Mullen-Hegland.\\
\end{center}
\end{strip}

Hereby  we can make the following preliminary conclusions.

\begin{enumerate}

\item The variable projection method requires a good initial guess for $\nu'$, as well as the mesh for $\nu$ must be fine enough. It is to be noted that for the  a real experiments it is technically difficult to employ the fine mesh. It is desirable to use derivatives  (Jacobi matrix), but this makes the method complicated. Finally, it is not easy to {\it a priory} estimate the number $n$ of spectral lines. But the false lines are slightly suppressed.

\item The algorithm of integral approximation does not use an initial guess for $\nu'$, but it also requires the fine mesh (for $\nu',$ which is easy to generate). The question of the true number of spectral lines is solved reliably, since the false line are safely filtered  (see. Fig. 7).

\end{enumerate}

In general, the two methods complement each other. The main difference is that  integral approximation is entirely linear method while  variable projection is a linear-nonlinear method.

%fig 7

\subsection{2nd example: real spectrum}

It is to be noted that in this example wavelength $\lambda$ is used instread of frequency
$\nu.$
Fig. 8 demonstrates the noisy spectrum of $\tilde{u}(\lambda)$ (20\% 
concentration of
carbon monoxide  CO (IRT themes) heated to a temperature of
$700\,^{\circ}{\rm C}$). The dataset is provided by Dept. Chem. and Biochem. Engg of DTU.
In this example we have the Dispersion type IF\cite{l9} (cf. Eq.~\ref{eqDT}) as follows
\begin{equation}
K(\lambda, \lambda') = \frac{a(\lambda)/2\pi}{(\lambda - \lambda')^2+[a(\lambda)/2]^2},
\label{eq32}
\end{equation}
where $a(\lambda) = q\lambda$ is the width of the instrumental function
 to half-power points, $q=0.0005$ such as $a(4780 \,\text{nm}) = 2.4\, \text{nm}.$

The following integral equation (cf. Eq.~\ref{eq2}) is solved
\begin{equation}
Az\equiv \int_a^b K(\lambda,\lambda') z(\lambda') \, d \lambda', \, c \leq \lambda \leq d,
\label{eq33}
\end{equation}
 the desired function $z(\lambda')$ is found using the quadrature method combined
with the regularization method (given by Eqs. \ref{eq3} -- \ref{eq6}) for 
$a = 4750,$ $b=4810,$ $c=4748,$ $d = 4812,$  $N=301,$
$m=321,$ $h=0.2,$ $ \alpha = \alpha_d = 10^{-2}.$  
As solution to SLAE Eq.~\ref{eq6} the regularized spectrum ${   z}_{\alpha}(\lambda')$ is
obtained (curve ``{\it 2}'' in Fig. 8).

Wavelengths  $\tilde{\lambda}_j'$ which correspond 
to the highest peaks (for $L=15$) are selected from the solution $z(\lambda')$.

The SLAE (cf. Eq. \ref{eq21})
\begin{equation}
\sum_{j=1}^L K(\lambda_i,\tilde{\lambda}_j') \tilde{z}_j +\tilde{F} = \tilde{u}(\lambda_i),
i=\overline{1,m}, \, c\leq \lambda_i \leq d
\label{eq34}
\end{equation} 
is solved, making it possible to refine amplitudes $\tilde{z}_j.$
The result is more accurate, see for example the vertical lines ``3'' in
Fig. 8. 

\begin{figure}[htbp]
\includegraphics[scale=0.467]{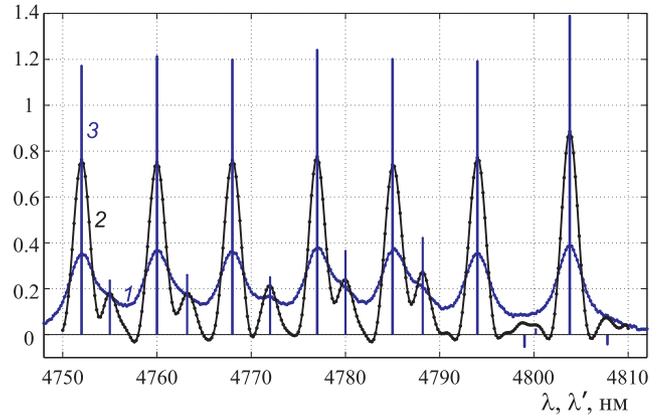}
\caption{Real  discrete spectrum. ``{\it 1}'' is measured noisy spectrum $\tilde{   u}(\lambda)$;
``{\it 2}'' is the spectrum ${   z}_{\alpha}(\lambda')$ obtained by the regularized quadrature method;
``{\it 3}'' is the reconstructed  discrete spectrum ${   z}_{j}(\tilde{\lambda}_j').$ }
\end{figure}

Based on the experience gained from this approach to solution of the problem, 
it can be concluded that 
7  largest lines have been safely recovered as well as 5 further weaker lines 
which have not been resolved in the measured spectrum $\tilde{u}(\lambda)$
 (see curve 1 in Fig. 8). As for the three  weakest
peaks in the high frequency part of the spectrum $z_{\alpha}(\lambda')$,  two of them 
have  negative
values, but  in the case of the maximum at $\lambda = 4800$ it is difficult to judge whether it corresponds to  real (weak) lines or to lines generated by the Gibbs phenomena.

\section{Conclusion}

We considered the problem of  effective enhancement of spectrometer resolution  without 
lines width artificial reduction.
The experimentû with synthetic and real data have demonstrated the high efficacy of the  proposed method. A substantial increase in the effective resolution of the spectrometer was observed.
The method enables enhancement of the resolution of spectral analysis through secure detection of closely located  lines, detection of weak lines in noisy spectra, revealing the fine structure of reconstructed discrete spectra  using mathematical and computer methods for spectral data processing. 

The proposed algorithm enjoy affordable computational complexity, elapsed time was about 1 sec in both examples. It can be implemented as an embedded software system within a spectrometer in order to deliver higher effective resolution.
Furthermore, an inexpensive spectrometer with a wide IF can be used combined with customised stand-alone Matlab software in order to enhance the original spectrometer's  effective resolution.

As a final remark it may be noted that the proposed algorithm for solving this inverse spectroscopy problem is a generic approach for analysing discrete spectra.
It can be used to recover smoothed and noisy spectra in various applications including  spectroscopy of gases, liquids, plasma and space objects, LiDARs, NMR spectroscopy, spectroscopic analysis of molten metal in blast furnaces etc.

\begin{acks}

The authors are grateful 
to the journal's referees for their comments on an earlier draft, to
Alexander Fateev (Dept. Chem. and Biochem. Engineering of DTU) for spectra measurements
and to Paul Leahy (University College Cork) for  helpful discussions.

\end{acks}
\begin{funding}
This work was funded by Russian Foundation of Basic Research project No.~13-08-00442
and by the International science and technology
cooperation program of China, project 2015DFR70850,  NSFC Grant No. 61673398.

\end{funding}

\end{document}